\documentclass[12pt]{amsart}
\usepackage{amsfonts}
\usepackage{epsfig}
\usepackage{graphicx}
\usepackage{amsmath}
\usepackage{amssymb}
\input{amssym.def}
\newtheorem{theorem}{Theorem}

\theoremstyle{remark}

\input epsf
\begin{document}
\title[Loewner equation]{Solutions of the Loewner equation with combined driving functions}
\author[D.~Prokhorov, A.~Zakharov, A.~Zherdev]{D.~Prokhorov, A.~Zakharov and A.~Zherdev}

\address{D.~Prokhorov: Department of Mathematics and Mechanics, Saratov State University, Saratov
410012, Russia} \email{ProkhorovDV@info.sgu.ru}
\address{A.~Zakharov: Department of Mathematics and Mechanics, Saratov State University, Saratov
410012, Russia} \email{Zaharovam@info.sgu.ru}
\address{A.~Zherdev: Department of Mathematics and Mechanics, Saratov State University, Saratov
410012, Russia} \email{Jerdevandrey@gmail.com}

\begin{abstract}
The paper is devoted to the multiple chordal Loewner differential equation with different driving functions on two time intervals. We obtain exact implicit or explicit solutions to the Loewner equations with piecewise constant driving functions and with combined constant and square root driving functions. In both cases, there is an analytical and geometrical description of generated traces.
\end{abstract}
\maketitle

\section{Introduction}

The Loewner differential equations \cite{Loewner} play important roles in the geometric function theory of complex analysis. We will discuss  a half-plane version of the Loewner equation, see e.g., [2, Chapter 4], generating self-maps of the upper half-plane $\mathbb H=\{z\in\mathbb C:\text{Im}\,z>0\}$. Given a simple curve $\Gamma$ in $\mathbb H$, emanating from a point on $\mathbb R$, and for an appropriate continuous parametrization $\Gamma(t)$ of $\Gamma$, $0\leq t\leq T$, there exists a unique conformal map $g(\cdot,t)$ from $\mathbb H\setminus\Gamma[0,t]$ onto $\mathbb H$ that obeys the hydrodynamic normalization near infinity, $$g(z,t)=z+\frac{2t}{z}+O\left(\frac{1}{|z|^2}\right),\;\;\;z\to\infty.$$ In this case, there is a continuous driving function $\lambda:[0,T]\to\mathbb R$ such that $g$ solves the chordal Loewner differential equation

\begin{equation}
\frac{\partial g(z,t)}{\partial t}=\frac{2}{g(z,t)-\lambda(t)},\;\;\;g(z,0)=z,\;\;\;0\leq t\leq T,\;\;\;z\in\mathbb H\setminus\Gamma[0,T]. \label{Loe1}
\end{equation}
We say that $g$ generates $\Gamma$.

If $\Gamma$ is a finite union of simple curves, probably with common points, we need to use the multiple Loewner differential equation $$\frac{\partial g(z,t)}{\partial t}=\sum_{k=1}^n\frac{2\mu_k}{g(z,t)-\lambda_k(t)},\;\;\;g(z,0)=z,\;\;\;0\leq t\leq T,\;\;\;z\in\mathbb H\setminus\Gamma[0,T],$$
with (piecewise) continuous driving functions $\lambda_k:[0,T]\to\mathbb R$ and positive numbers $\mu_k$, $k=1,\dots,n$, $\sum_{k=1}^n\mu_k=1$.

In this paper, we restrict ourselves to $n=2$ and $\mu_1=\mu_2=\frac{1}{2}$. So we consider the Loewner differential equation
\begin{equation}
\frac{\partial g(z,t)}{\partial t}=\sum_{k=1}^2\frac{1}{g(z,t)-\lambda_k(t)},\;\;\;g(z,0)=z,\;\;\;0\leq t\leq T,\;\;\;z\in\mathbb H\setminus\Gamma[0,T]. \label{Loe2}
\end{equation}

There are only few known examples of driving functions in (\ref{Loe1}) or (\ref{Loe2}) that admit explicit integration of this equation and describe corresponding traces $\Gamma$. In \cite{KagNieKad}, the authors solve equation (\ref{Loe2}) with constant driving functions $\lambda_1<0$ and $\lambda_2=-\lambda_1$. For $n=1$, the authors of \cite{KagNieKad} give a full description of a trace generated in equation (\ref{Loe1}) if the driving function has the form $\lambda(t)=A\sqrt t$, $A>0$.

We aim to develop integration possibilities for equation (\ref{Loe2}) with combined driving functions $\lambda_1$ and $\lambda_2$, $\lambda_1=-\lambda_2$, when

\begin{equation}
\lambda_2(t) =
\begin{cases}\displaystyle
0, & 0\leq t<t_0, \\
A, & t_0\leq t\leq T,
\end{cases} \label{dri1}
\end{equation}
or

\begin{equation}
\lambda_2(t) =
\begin{cases}\displaystyle
0, & 0\leq t<t_0, \\
A\sqrt{t-t_0}, & t_0\leq t\leq T,
\end{cases} \label{dri2}
\end{equation}
for arbitrary $A>0$ and $t_0>0$ and a certain $T>t_0$.

Note that both driving functions $\lambda_1,\lambda_2$ in (\ref{dri2}) are continuous on $[0,T]$ while driving functions $\lambda_1,\lambda_2$ in (\ref{dri1}) have jumps at $t_0$.

In Section 2, we integrate Loewner equation (\ref{Loe2}) with piecewise constant driving functions (\ref{dri1}), see Theorem 1, and show that a solution $g(\cdot,t)$ maps $\mathbb H\setminus\Gamma$ onto $\mathbb H$ where $\Gamma$ is a union of a segment $\Gamma_0$ on the upper imaginary half-axis and a pair of curves $\Gamma_1$ and $\Gamma_2$ which are symmetric with respect to the imaginary axis and emanate either from points on $\mathbb R$ if $A>2\sqrt{t_0}$ or from points on $\Gamma_0$ if $A<2\sqrt{t_0}$. If $A=2\sqrt{t_0}$, the boundary symmetric curves $\Gamma_1$ and $\Gamma_2$ emanate from the origin under angles $\pm\frac{\pi}{4}$ to the real axis $\mathbb R$. We give implicit representations of $\Gamma_1$ and $\Gamma_2$ and asymptotic expansions for $\Gamma_1$ and $\Gamma_2$ near $t=t_0$.

In Section 3, we integrate Loewner equation (\ref{Loe2}) with continuous driving functions (\ref{dri2}) which are constant on $[0,t_0)$ and square root functions on $[t_0,T]$, see Theorem 2. We show that a solution $g(\cdot,t)$ maps $\mathbb H\setminus\Gamma$ onto $\mathbb H$ where $\Gamma$ is a union of the segment $[0,i2\sqrt{t_0}]$ and a pair of curves which are symmetric with respect to the imaginary axis and emanate from the point $i2\sqrt{t_0}$. We give explicit representations of boundary curves and their asymptotic expansions near $t=t_0$.

In Section 4, we discuss an interrelation between exact solutions for the standard Loewner equation on two separate time intervals and its multiple version.

\section{Loewner equation with piecewise constant driving functions}

Let us solve the multiple Loewner differential equation (\ref{Loe2}) with combined driving functions (\ref{dri1}) that are piecewise constant on $[0,T]$.

\begin{theorem}
There exists $T>t_0$ for which the multiple Loewner differential equation (\ref{Loe2}) with combined driving functions (\ref{dri1}) has a solution $w=g(z,t)$ on $[0,T]$. On $[0,t_0]$, $g(z,t)=\sqrt{z^2+4t}$, and on $[t_0,T]$, $w=g(z,t)$ satisfies the implicit equation
\begin{equation}
w^2-z^2-A^2\log\frac{w^2}{z^2+4t_0}=4t,\,\,\,g(z,t_0)=\sqrt{z^2+4t_0},\;\;\;z\in\mathbb H\setminus[0,i2\sqrt{t_0}], \label{Gam1}
\end{equation}
where the continuous branches of $\log w$ and $\log z$ are real when $w$ and $z$ are positive. The function $g(z,T)$ maps $\mathbb H\setminus\Gamma$ onto $\mathbb H$ according to the following three cases:

(i) If $A>2\sqrt{t_0}$, then $\Gamma=\cup_{k=0}^2\Gamma_k$, where $\Gamma_0=[0,2i\sqrt{t_0}]$, $\Gamma_2[0,T]$ is a curve which emanates from $\sqrt{A^2-4t_0}$ and is orthogonal to $\mathbb R$ at this point, $\Gamma_1[0,T]$ is symmetric to $\Gamma_2[0,T]$ with respect to the imaginary axis;

(ii) If $A<2\sqrt{t_0}$, then $\Gamma=\cup_{k=0}^2\Gamma_k^*$, where $\Gamma_0^*=\Gamma_0$, $\Gamma_2^*[0,T]$ is a curve which emanates from $i\sqrt{4t_0-A^2}$ and is orthogonal to the imaginary axis at this point, $\Gamma_1^*[0,T]$ is symmetric to $\Gamma_2^*[0,T]$ with respect to the imaginary axis;

(iii) If $A=2\sqrt{t_0}$, then $\Gamma=\cup_{k=0}^2\Gamma_k^{**}$, where $\Gamma_0^{**}=\Gamma_0$, $\Gamma_2^{**}[0,T]$ is a curve which emanates from the origin under the angle $\frac{\pi}{4}$ to $\mathbb R$, $\Gamma_1^{**}[0,T]$ is symmetric to $\Gamma_2^{**}[0,T]$ with respect to the imaginary axis.

\end{theorem}

\begin{figure}
\epsfig{file=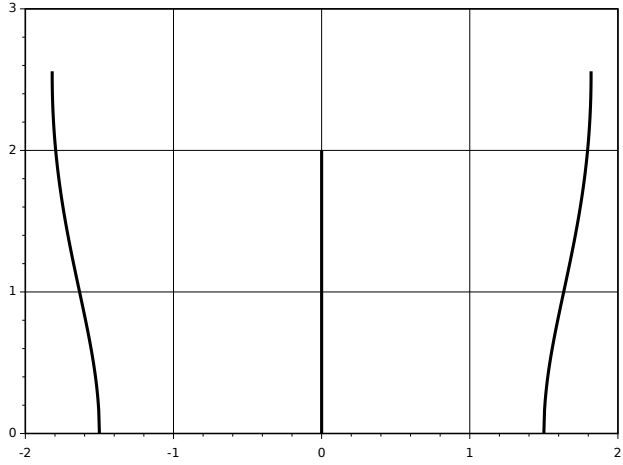, width=3.2in}
\caption{$\Gamma$ for $t_0=1, T=3, A=2.5 \,(i)$}
\end{figure}
\begin{figure}
\epsfig{file=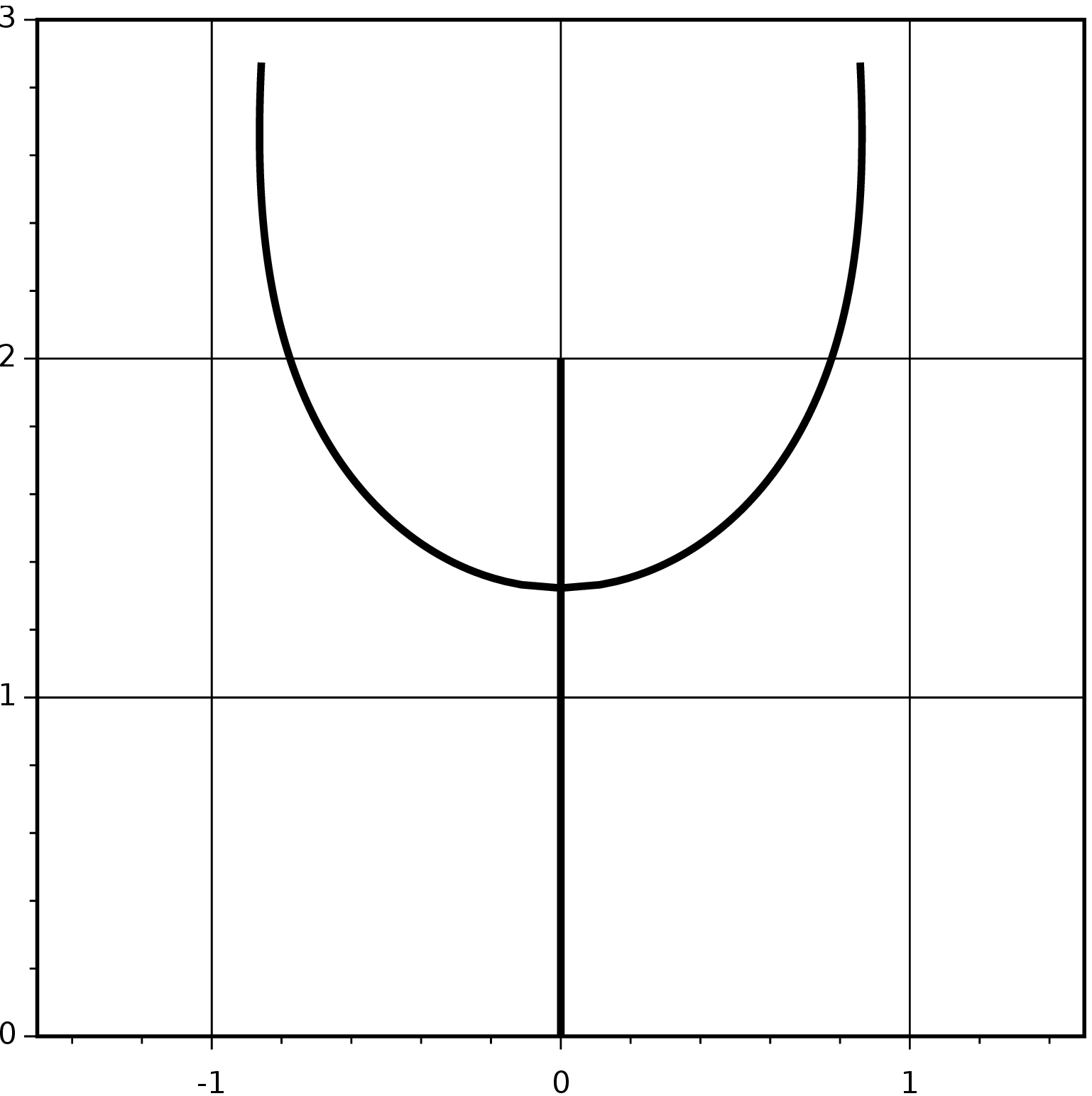, width=2.4in}
\caption{$\Gamma$ for $t_0=1, T=3, A=1.5 \,(ii)$}
\end{figure}
\begin{figure}
\epsfig{file=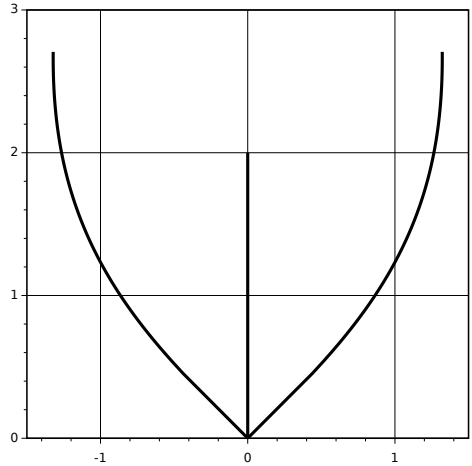, width=2.4in}
\caption{$\Gamma$ for $t_0=1, T=3, A=2 \,(iii)$}
\end{figure}

\begin{proof}
It is a well-known result on $[0,t_0]$ that $g(z,t)=\sqrt{z^2+4t}$, $z\in\mathbb H$, see, e.g., [2, p.95], \cite{KagNieKad}. Next, we have to solve the multiple Loewner equation
\begin{equation}
\frac{dw}{dt}=\frac{1}{w+A}+\frac{1}{w-A}=\frac{2w}{w^2-A^2},\;\;\;w(z,t_0)=g(z,t_0),\;\;\;t_0\leq t\leq T. \label{con1}
\end{equation}
The function $g(z,t_0)$ maps $\mathbb H\setminus[0,i2\sqrt{t_0}]$ onto $\mathbb H$. Differential equation (\ref{con1}) with separated variables $w$ and $z$ has a general solution in the form $$w^2-2A^2\log w=4t+c$$ with an arbitrary constant $c$. The initial value allows us to determine $c$ as $$c=z^2-A^2\log(z^2+4t_0).$$ So we find an implicit solution $w=w(z,t)$ to the Cauchy problem (\ref{con1}) as it is presented in (\ref{Gam1}).

Differential equation (\ref{con1}) generates two traces $\Gamma_1$ and $\Gamma_2$ that are symmetric with respect to the imaginary axis and emanate from two points $g(z,t_0)=\pm A$ on $\mathbb R$. Let $\Gamma_2$ correspond to $g(z,t_0)=A$ and let $\Gamma_2$ be given by $z=z(t)$. Then the line of singularities $z(t)$ satisfies the equation $$w(z(t),t)=A,\;\;\;t\geq t_0.$$ Together with (\ref{Gam1}) this leads to the equality
\begin{equation}
A^2-z^2(t)-A^2\log\frac{A^2}{z^2(t)+4t_0}=4t,\;\;\;t\geq t_0,\;\;\;g(z(t_0),t_0)=A. \label{con2}
\end{equation}
The equality $g(z(t_0),t_0)=A$ is equivalent to $z(t_0)=\sqrt{A^2-4t_0}$. A disposition of the initial point of $\Gamma_2$ depends on the sign of $A^2-4t_0$. Consider three possible cases.

{\bf Case (i):} $A>2\sqrt{t_0}$. Then $z(t_0)>0$ and $\Gamma_2$ emanates from the point on the positive real half-axis.

{\bf Case (ii):} $0<A<2\sqrt{t_0}$. Then $z(t_0)$ is pure imaginary and $\Gamma_2$ emanates from the point on $(0,i2\sqrt{t_0})$.

{\bf Case (iii):} $A=2\sqrt{t_0}$. Then $z(t_0)=0$ and $\Gamma_2$ emanates from the origin.

Equality (\ref{con2}) is an implicit representation of $\Gamma_2$. Find an asymptotic expansion of $\Gamma_2$ near the initial point in all the three cases.

Differentiate (\ref{con2}) and obtain
\begin{equation}
(z^2(t))'=\frac{4(z^2(t)+4t_0)}{A^2-4t_0-z^2(t)},\;\;\;z^2(t_0)=A^2-4t_0,\;\;\;t\geq t_0. \label{con3}
\end{equation}

This allows us to find an asymptotic expansion for $z(t)$ near $t_0$. In cases (i) and (ii) it is reasonable to set $$z(t)=\sqrt{A^2-4t_0}+a\sqrt{t-t_0}+o(\sqrt{t-t_0}),\;\;\;t\to t_0^+.$$ Hence
$$(z^2(t))'=\frac{a\sqrt{A^2-4t_0}}{\sqrt{t-t_0}}+o\left(\frac{1}{\sqrt{t-t_0}}\right),\;\;\;t\to t_0^+.$$
Substitute expansions for $z(t)$ and $(z^2(t))'$ in (\ref{con3}) and see that $$\frac{a\sqrt{A^2-4t_0}}{\sqrt{t-t_0}}=-\frac{2A^2}{a\sqrt{A^2-4t_0}\sqrt{t-t_0}}+o\left(\frac{1}{\sqrt{t-t_0}}\right),\;\;\;t\to t_0^+,$$ which gives that $$a^2=-\frac{2A^2}{A^2-4t_0}.$$

In Case (i) $A^2>4t_0$: $$z(t)=\sqrt{A^2-4t_0}+i\frac{\sqrt2A}{\sqrt{A^2-4t_0}}\sqrt{t-t_0}+o(\sqrt{t-t_0}),\;\;\;t\to t_0^+.$$ So $z(t)$ is orthogonal to $\mathbb R$ at $z=\sqrt{A^2-4t_0}$.

In Case (ii) $A^2<4t_0$: $$z(t)=i\sqrt{4t_0-A^2}+\frac{\sqrt2A}{\sqrt{4t_0-A^2}}\sqrt{t-t_0}+o(\sqrt{t-t_0}),\;\;\;t\to t_0^+.$$ So $z(t)$ is orthogonal to the imaginary axis at $z=i\sqrt{4t_0-A^2}$.

Case (iii) $A^2=4t_0$ requires another asymptotic behavior of the trace near the origin. Formula (\ref{con3}) transforms to the following $$(z^4(t))'=-8(z^2(t)+4t_0),\;\;\;z(t_0)=0.$$ Present another reasonable asymptotic expansion for $z(t)$, $$z(t)=b\root4\of{t-t_0}+o(\root4\of{t-t_0}),\;\;\;t\to t_0^+.$$ Take into account both last formulas and obtain that $b^4=-32t_0.$ Thus $$z(t)=e^{i\frac{\pi}{4}}2\root4\of{2t_0}\root4\of{t-t_0}+o(\root4\of{t-t_0}),\;\;\;t\to t_0^+.$$ So $z(t)$ is tangential to the radial ray under the angle $\frac{\pi}{4}$ to $\mathbb R$ at the origin.

Similarly to (\ref{con2}) derive an implicit representation for $z(t)$ in Case (iii). Integrate the differential equation for $(z^4(t))'$ to get the needed equation $$z^2+4\log\frac{A^2}{z^2+4t_0}=A^2-4t,\;\;\;z(t_0)=0.$$

The boundary curve $\Gamma_1$ can be studied similarly. However, it is symmetric to $\Gamma_1$ with respect to the imaginary axis due to the symmetric disposition of points $\pm A$ and the symmetric trace on the time segment $[0,t_0]$.

It remains to observe what does occur with the boundary $[0,i2\sqrt{t_0}]$ when $t$ varies along $[t_0,T]$. The implicit representation (\ref{Gam1}) implies that the two singular points $w=0$ and $z=i2\sqrt{t_0}$ appear simultaneously. As far as $w=0$ is constant on $[t_0,T]$ according to (\ref{con1}), the corresponding $z=i2\sqrt{t_0}$ also does not move for $t$ on $[t_0,T]$. Inner points of the segment $[0,i2\sqrt{t_0}]$ cannot leave the imaginary axis because of symmetric properties of conformal mappings generated by symmetric driving functions (\ref{dri1}). This means that the segment $[0,i2\sqrt{t_0}]$ is a part of the boundary set $\Gamma$, and there are no additional parts of $\Gamma$ on the imaginary axis, which completes the proof of Theorem 1.
\end{proof}

\section{Combined constant and square root driving functions}

Now we will solve the multiple Loewner differential equation (\ref{Loe2}) with combined driving functions (\ref{dri2}) that are continuous on $[0,T]$.

\begin{theorem} For every $T>t_0$, the multiple Loewner differential equation (\ref{Loe2}) with combined driving functions (\ref{dri2}) has a solution $w=g(z,t)$ on $[0,T]$. On $[0,t_0]$, $g(z,t)=\sqrt{z^2+4t}$, and on $[t_0,T]$, $w=g(z,t)$ satisfies the implicit equation $$(A^2+4)(t-t_0)=w^2-(z^2+4t_0)^{\frac{A^2}{4}+1}w^{-\frac{A^2}{2}},\;\;\;w(z,t_0)=\sqrt{z^2+4t_0},$$ where the branches of power functions are such that they are positive when $z^2+4t_0$ and $w$ are positive. The function $g(z,T)$ maps $\mathbb H\setminus\Gamma$ onto $\mathbb H$, $\Gamma=\cup_{k=0}^2\Gamma_k$, where $\Gamma_0$ is the segment $[0,i2\sqrt{t_0}]$, $\Gamma_2$ is a square root of a rectilinear segment under the angle $4\pi/(A^2+4)$ to $\mathbb R$ from $(-4t_0)$, and $\Gamma_1$ is symmetric to $\Gamma_2$ with respect to the imaginary axis.
\end{theorem}

\begin{figure}
\epsfig{file=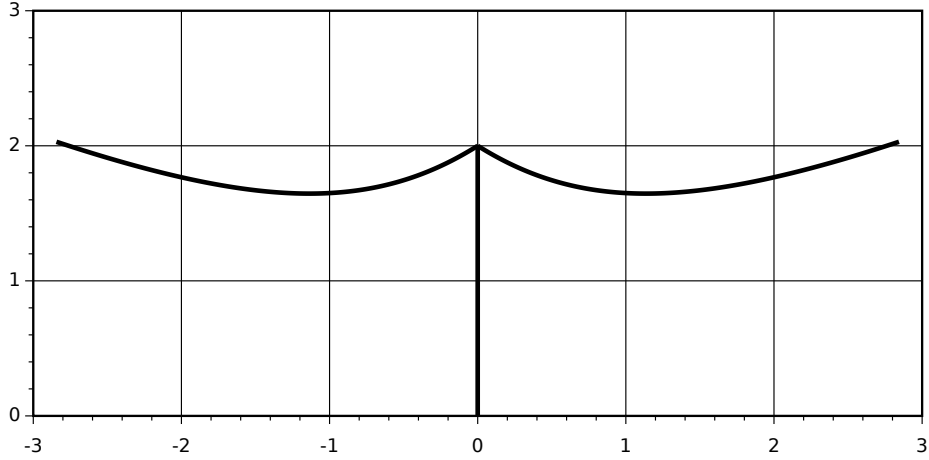, width=4.8in}
\caption{$\Gamma$ for $t_0=1, T=3, A=3.0$}
\end{figure}

\begin{proof} As in Theorem 1, on $[0,t_0]$, the function $g(z,t)=\sqrt{z^2+4t}$, $z\in\mathbb H$, solves the chordal Loewner differential equation (\ref{Loe2}) with vanishing driving functions. Next, we have to solve the multiple Loewner equation
\begin{equation}
\frac{dw}{dt}=\frac{1}{w+A\sqrt{t-t_0}}+\frac{1}{w-A\sqrt{t-t_0}}=\frac{2w}{w^2-A^2(t-t_0)},\;\;\;t_0\leq t\leq T, \label{con4}
\end{equation}
with the initial condition $w(z,t_0)=g(z,t_0)$. Remind that $g(z,t_0)=\sqrt{z^2+4t_0}$ maps $\mathbb H\setminus[0,2i\sqrt{t_0}]$ onto $\mathbb H$.

Note that differential equation (\ref{con4}) is linear with respect to $t$. Its general solution is given implicitly by $$t-t_0=\frac{w^2}{A^2+4}+cw^{-\frac{A^2}{2}}$$ with an arbitrary constant $c$. The initial value allows us to determine $c$ from the equation $$0=\frac{z^2+4t_0}{A^2+4}+c(z^2+4t_0)^{-\frac{A^2}{4}}$$ so that $$c=-\frac{(z^2+4t_0)^{\frac{A^2}{4}+1}}{A^2+4}.$$ So we find an implicit solution to the Cauchy problem (\ref{con4}) as
\begin{equation}
t-t_0=\frac{w^2}{A^2+4}-\frac{(z^2+4t_0)^{\frac{A^2}{4}+1}}{A^2+4}w^{-\frac{A^2}{2}}, \label{con5}
\end{equation}
which proves the first statement of Theorem 2 for a certain $T>t_0$.

Differential equation (\ref{con4}) generates two traces $\Gamma_1$ and $\Gamma_2$ symmetric with respect to the imaginary axis and emanating from the common point $g(0,t_0)=i2\sqrt{t_0}$ on the imaginary axis. Let $\Gamma_2$ be situated in the right half-plane for a certain $T>t_0$ and let $\Gamma_2$ be given by $z(t)$. Then the line of singularities $z(t)$ satisfies the equation $$w(z(t),t)=A\sqrt{t-t_0},\;\;\;t\geq t_0.$$ Together with (\ref{con5}) this leads to the equality $$t-t_0=\frac{A^2(t-t_0)}{A^2+4}-\frac{(z^2(t)+4t_0)^{\frac{A^2}{4}+1}}{A^2+4}(A\sqrt{t-t_0})^{-\frac{A^2}{2}},\;\;\;t\geq t_0.$$ Transform this expression to the following $$4(t-t_0)^{\frac{A^2}{4}+1}=-(z^2(t)+4t_0)^{\frac{A^2}{4}+1}A^{-\frac{A^2}{2}}$$ and give the explicit formula for $z(t)$, $$z(t)=\left[e^{\frac{i4\pi}{A^2+4}}2^{\frac{8}{A^2+4}}A^{\frac{2A^2}{A^2+4}}(t-t_0)-4t_0\right]^{\frac{1}{2}},\;\;\;t\geq t_0.$$

It is worth noting that $z(t)$ is the square root of a rectilinear segment under the angle $4\pi/(A^2+4)$ to $\mathbb R$ from $(-4t_0)$. The slope of the rectilinear segment is changing from $\pi$ to 0 when $A$ is growing from 0 to infinity. Therefore, $\Gamma_2$ emanates from $i2\sqrt{t_0}$ and is tangential to the ray under the angle $$\frac{\pi 4}{A^2+4}-\frac{\pi}{2}=\frac{\pi(4-A^2)}{2(A^2+4)}$$ to $\mathbb R$ at this endpoint. The slope of $\Gamma_2$ is changing from $\pi/2$ to $(-\pi/2)$ when $A$ is growing from 0 to infinity.

The last reasoning explains that $\Gamma_2$ stays in the right half-plane for all $A$ and $T>t_0$ and it is a simple curve.

The boundary curve $\Gamma_1$ can be studied similarly. However, it is symmetric to $\Gamma_1$ with respect to the imaginary axis due to the symmetric properties of the driving functions $\pm A\sqrt{t-t_0}$ and the symmetric trace on the time segment $[0,t_0]$.

It is known that the Loewner differential equation generates simple traces up to the moment $t$ when either lines of singularities $\Gamma$ meet the real axis $\mathbb R$ or $\Gamma$ has self-intersection, see, e.g., \cite{Lind}. We showed that, under conditions of Theorem 2, the curve $\Gamma_2$ stays in the right half-plane for all $t$ and does not reach $\mathbb R$. Similarly, the curve $\Gamma_1$ stays in the left half-plane and does not reach $\mathbb R$. Both $\Gamma_2$ and $\Gamma_1$ do not meet $\Gamma_0:=[0,i2\sqrt{t_0}]$. Hence the Loewner generating process develops in time for all $T>t_0$. This completes the proof of Theorem 2.

\end{proof}

\section{Conclusions}

The proofs of Theorems 1 and 2 are based on the knowledge of integrability cases of the Loewner differential equation for constant and square root driving functions both in the standard and multiple versions. There are some more known driving functions that admit explicit or implicit integration of the Loewner equation. Therefore, it is possible to present new examples of combined driving functions in the Loewner equation with several contact points which join different driving functions and lead to exact solutions.

Point out at such examples. Besides constant and square root driving functions, Kager, Nienhuis and Kadanoff \cite{KagNieKad} considered linear driving functions and obtained exact solutions of the Loewner equation. We have to add that its adaptation to the multiple equation is not so successful in getting exact solutions.

In \cite{ProZak}, the authors found an implicit exact solution of the Loewner equation with the exponential driving function $A(e^t-1)$. Moreover, this driving function is well-adapted to express an exact solution for the multiple Loewner equation.

There is another approach in the exact solution problem when driving functions are determined for given traces of the Loewner equation. We refer to \cite{ProVas} where the problem was solved for the circular arc in $\mathbb H$ tangential to $\mathbb R$ at 0. This result was generalized in \cite{LauWu} for powers of this arc and in \cite{WuJiaDon} for tangential curves close to this arc. It was proved in \cite{ProVas} that the tangential circular arc of radius 1 and centered at $i$ is generated by the driving function $\lambda(t)=3\alpha(t)+2\sqrt{-\alpha(t)\pi}$ where $\alpha=\alpha(t)$ is an algebraic function satisfying the equation $$\alpha(3\alpha+4\sqrt{-\alpha\pi})=-6t,\;\;\;t\geq0.$$ A similar problem was solved by Wu in \cite{Wu} for circular arcs in $\mathbb H$ which meet $\mathbb R$ orthogonally.

\end{document}